\documentclass[11pt,a4paper,twoside]{article}
\usepackage{amssymb}
\usepackage{amsmath}
\usepackage{fancyhdr}
\usepackage{amsthm}
\usepackage{rotate}
\usepackage{graphicx}
\usepackage[T1]{fontenc}
\usepackage{afterpage}
\theoremstyle{definition}
\newtheorem{definition}{Definition}[section]
\newtheorem{example}[definition]{Example}
\theoremstyle{plain}
\newtheorem{theorem}[definition]{Theorem}
\newtheorem{lemma}[definition]{Lemma}
\newtheorem{corollary}[definition]{Corollary}
\newtheorem{proposition}[definition]{Proposition}
\newtheorem{problem}[definition]{Problem}
\newtheorem{conjecture}[definition]{Conjecture}
\newtheorem{remark}[definition]{Remark}

\theoremstyle{remark}

\newcommand{\octo}[1]{\mathbb O(#1)}
\newcommand{\vm}[4]{\left(\begin{array}{cc}
    #1&#2\\#3&#4\end{array}\right)}
\newcommand{\ol}[1]{\overline{#1}}
\newcommand{\paige}[1]{M^*(#1)}
\newcommand{\id}{\mathrm{id}}
\DeclareMathOperator{\chr}{char}
\DeclareMathOperator{\discr}{discr}
\DeclareMathOperator{\Inn}{Inn}
 \DeclareMathOperator{\Mlt}{Mlt}
\DeclareMathOperator{\Coll}{Coll}
 \DeclareMathOperator{\Aut}{Aut}
\begin{document}

\afterpage{\rhead[]{\thepage} \chead[\small G.~P.~Nagy and
P.~Vojt\v{e}chovsk\'y]{\small  Octonions and simple Moufang loops}
\lhead[\thepage]{} }
\begin{center}
\vspace*{2pt}
{\Large \textbf{Octonions, simple Moufang loops and triality}}\\[36pt]
{\large \textsf{\emph{G\'abor P. Nagy and Petr
Vojt\v{e}chovsk\'y}}}
\\[36pt]

\textbf{Abstract}
\end{center}
{\footnotesize Nonassociative finite simple Moufang loops are
exactly the loops constructed by Paige from Zorn vector matrix
algebras. We prove this result anew, using geometric loop theory.
In order to make the paper accessible to a broader audience, we
carefully discuss the connections between composition algebras,
simple Moufang loops, simple Moufang $3$-nets, $S$-simple groups
and groups with triality. Related results on multiplication
groups, automorphisms groups and generators of Paige loops are
provided.}

 \footnote{

$\rule{3mm}{0mm}$group, group with triality, net }\\
\footnote{$\rule{-6mm}{0mm}$The first author was supported by the
FKFP grant $0063/2001$ of the Hungarian Ministry for Education and
the OTKA grants nos. F042959 and T043758. The second author
partially supported by the Grant Agency of Charles University,
grant number $269/2001/$B-MAT/MFF.}

\hyphenation{col-lin-ea-tion col-lin-ea-tions} {\small
\tableofcontents }
\section{Introduction}

The goal of this paper is to present the classification of finite
simple Moufang loops in an accessible and uniform way to a broad
audience of researchers in nonassociative algebra. The results are
not new but the arguments often are. Although not all proofs are
included, our intention was to leave out only those proofs that
are standard (that is those that can be found in many sources),
those that are purely group-theoretical, and those that require
only basic knowledge of loop theory. We have rewritten many proofs
using geometric loop theory---a more suitable setting for this
kind of reasoning. To emphasize the links to other areas of loop
theory and algebra, we comment on definitions and results
generously, although most of the remarks we make are not essential
later in the text.

Here is a brief description of the content of this paper. After
reviewing some basic properties of loops, nets and composition
algebras, we construct a family of simple Moufang loops from the
Zorn alternative algebras. These loops are also known as Paige
loops. We then briefly discuss the multiplication groups of Paige
loops, because these are essential in the classification.

With every Moufang loop we associate a Moufang $3$-net, and with
this $3$-net we associate a group with triality. An
$S$-homomorphism is a homomorphism between two groups with
triality that preserves the respective triality automorphisms.
This leads us to the concept of $S$-simple groups with triality,
which we classify. The group with triality $G$ associated with a
simple Moufang loop $L$ must be $S$-simple. Moreover, when $L$ is
nonassociative $G$ must be simple. This is the moment when we use
results of Liebeck concerning the classification of finite simple
groups with triality. His work is based on the classification of
finite simple groups. The fact that there are no other
nonassociative finite simple Moufang loops besides finite Paige
loops then follows easily.

Building on the geometric understanding we have obtained so far,
we determine the automorphism groups of all Paige loops
constructed over perfect fields. We conclude the paper with
several results concerning the generators of finite Paige loops
and integral Cayley numbers. All these results are mentioned
because they point once again towards classical groups. Several
problems and conjectures are pondered in the last section.

A few words concerning the notation: As is the habit among many
loop theorists, we write maps to the right of their arguments, and
therefore compose maps from left to right. The only exception to
this rule are some traditional maps, such as the determinant
$\det$. A subloop generated by $S$ will be denoted by $\langle
S\rangle$. The symmetric group on $n$ points is denoted by $S_n$.

\section{Loops and nets}\label{Sc:LN}

\noindent We now give a brief overview of definitions and results
concerning loops and nets. Nets (also called webs in the
literature) form the foundations of the geometric loop theory. All
material covered in \ref{Ss:LN1}--\ref{Ss:LN3} can be found in
\cite{Bruck} and \cite{Pflug}, with proofs. We refer the reader to
\cite[Ch.\ II\,]{Pflug} and \cite[Ch.\ VIII,\,X\,]{CPS} for
further study of nets.

\subsection{Quasigroups and loops}\label{Ss:LN1}

Let $Q=(Q,\cdot)$ be a groupoid. Then $Q$ is a \emph{quasigroup}
if the equation $x\cdot y = z$ has a unique solution in $Q$
whenever two of the three elements $x$, $y$, $z\in Q$ are
specified. Quasigroups are interesting in their own right, but
also appear in combinatorics under the name \emph{latin squares}
(more precisely, multiplication tables of finite quasigroups are
exactly latin squares), and in universal algebra, where
subvarieties of quasigroups are often used to provide an instance
of some universal algebraic notion that cannot be demonstrated in
groups or other rigid objects. We ought to point out that in order
to define the variety of quasigroups equationally, one must
introduce additional operations $\backslash$ and $/$ for left and
right division, respectively.

A quasigroup $Q$ that possesses an element $e$ satisfying $e\cdot
x = x\cdot e = x$ for every $x\in Q$ is called a \emph{loop} with
\emph{neutral element} $e$. The vastness of the variety of loops
dictates to focus on some subvariety, usually defined by an
identity approximating the associative law. (Associative loops are
exactly groups.) In this paper, we will be concerned with
\emph{Moufang loops}, which are loops satisfying any one of the
three equivalent \emph{Moufang identities}
\begin{equation}\label{Eq:Moufang}
    ((xy)x)z = x(y(xz)), \quad
    ((xy)z)y = x(y(zy)), \quad
    (xy)(zx) = (x(yz))x,
\end{equation}
and in particular with simple Moufang loops (see below). Every
element $x$ of a Moufang loop is accompanied by its
\emph{two-sided inverse} $x^{-1}$ satisfying $xx^{-1}=x^{-1}x=e$.
Any two elements of a Moufang loop generate a subgroup, and thus
$(xy)^{-1}=y^{-1}x^{-1}$.

Each element $x$ of a loop $Q$ gives rise to two permutations on
$Q$, the \emph{left translation} $L_x:y\mapsto xy$ and the
\emph{right translation} $R_x:y\mapsto yx$. The group $\Mlt{Q}$
generated by all left and right translations is known as the
\emph{multiplication group} of $Q$. The subloop $\Inn{Q}$ of
$\Mlt{Q}$ generated by all maps $L_xL_yL_{yx}^{-1}$,
$R_xR_yR_{xy}^{-1}$ and $R_xL_x^{-1}$, for $x$, $y\in Q$, is
called the \emph{inner mapping group} of $Q$. It consists of all
$\varphi\in\Mlt{Q}$ such that $e\varphi=e$.

A subloop $S$ of $Q$ is \emph{normal} in $Q$ if $S\varphi = S$ for
every $\varphi\in\Inn{Q}$. The loop $Q$ is said to be
\emph{simple} if the only normal subloops of $Q$ are $Q$ and
$\{e\}$.

In any loop $Q$, the \emph{commutator} of $x$, $y\in Q$ is the
unique element $[x,y]\in Q$ satisfying $xy=(yx)[x,y]$, and the
\emph{associator} of $x$, $y$, $z\in Q$ is the unique element
$[x,y,z]\in Q$ satisfying $(xy)z = (x(yz))[x,y,z]$. We prefer to
call the subloop $C(Q)$ of $Q$ consisting of all elements $x$ such
that $[x,y]=[y,x]=e$ for every $y\in Q$ the \emph{commutant} of
$Q$. (Some authors use the name \emph{centrum} or \emph{Moufang
center}.) The subloop $N(Q)$ consisting of all $x\in Q$ such that
$[x,y,z]=[y,x,z]=[y,z,x]=e$ holds for every $y$, $z\in Q$ is known
as the \emph{nucleus} of $Q$. Then $Z(Q)=C(Q)\cap N(Q)$ is the
\emph{center} of $Q$, which is always a normal subloop of $Q$.

\subsection{Isotopisms versus isomorphisms}\label{Ss:LN2}

Quasigroups and loops can be classified up to isomorphism or up to
isotopism. When $Q_1$, $Q_2$ are quasigroups, then the triple
$(\alpha$, $\beta$, $\gamma)$ of bijections from $Q_1$ onto $Q_2$
is an \emph{isotopism} of $Q_1$ onto $Q_2$ if $x\alpha\cdot y\beta
= (x\cdot y)\gamma$ holds for every $x$, $y\in Q_1$. An isotopism
with $Q_1=Q_2$ is called an \emph{autotopism}. Every isomorphism
$\alpha$ gives rise to an isotopism $(\alpha$, $\alpha$,
$\alpha)$. The notion of isotopism is superfluous in group theory,
as any two groups that are isotopic are already isomorphic.

In terms of multiplication tables, $Q_1$ and $Q_2$ are isotopic if
the multiplication table of $Q_2$ can be obtained from the
multiplication table of $Q_1$ by permuting the rows (by $\alpha$),
the columns (by $\beta$), and by renaming the elements (by
$\gamma$). Isotopisms are therefore appropriate morphisms for the
study of quasigroups and loops. On the other hand, every
quasigroup is isotopic to a loop, which shows that the algebraic
properties of isotopic quasigroups can differ substantially.
Fortunately, the classification of finite simple Moufang loops is
the same no matter which kind of equivalence (isotopism or
isomorphism) we use. This is because (as we shall see) there is at
most one nonassociative finite simple Moufang loop of a given
order, up to isomorphism.

A loop $L$ is a \emph{$G$-loop} if every loop isotopic to $L$ is
isomorphic to $L$. So, finite simple Moufang loops are $G$-loops.

\subsection{Loops and $3$-nets}\label{Ss:LN3}

Let $k>2$ be an integer, $\mathcal P$ a set, and $\mathcal L_1,
\ldots, \mathcal L_k$ disjoint sets of subsets of $\mathcal P$.
Put $\mathcal L = \bigcup \mathcal L_i$. We call the elements of
$\mathcal P$ and $\mathcal L$ {\em points} and {\em lines},
respectively, and use the common geometric terminology, such as
``all lines through the point $P$'', etc. For $\ell \in \mathcal
L_i$, we also speak of a \emph{line of type $i$} or an
\emph{$i$-line}. Lines of the same type are called {\em parallel}.

The pair $(\mathcal P, \mathcal L)$ is a {\em $k$-net} if the
following axioms hold:\\
$\rule{3mm}{0mm}1)\rule{0mm}{6mm}$ Distinct lines of the same type
are disjoint.\\
$\rule{3mm}{0mm}2)\rule{0mm}{6mm}$ Two lines of different types
have precisely one point in common.\\
$\rule{3mm}{0mm}3)\rule{0mm}{6mm}$ Through any point, there is
precisely one line of each type.\\

Upon interchanging the roles of points and lines, we obtain {\em
dual $k$-nets}. In that case, the points can be partitioned into
$k$ classes so that:\\
$\rule{3mm}{0mm}1')\rule{0mm}{6mm}$ Distinct points of the same
type are not connected by a line.\\
$\rule{3mm}{0mm}2')\rule{0mm}{6mm}$ Two points of different types
are connected by a unique line.\\
$\rule{3mm}{0mm}3')\rule{0mm}{6mm}$ Every line consists of $k$
points of pairwise different types.\\

There is a natural relation between loops and $3$-nets. Let us
first start from a loop $L$ and put $\mathcal P=L \times L$.
Define the line classes

\[\begin{array}{rcl}
{\mathcal L}_1&=&\{ \{(x,c) \mid x\in L\} \mid c\in L\}, \\
{\mathcal L}_2&=&\{ \{(c,y) \mid y\in L\} \mid c\in L\}, \\
{\mathcal L}_3&=&\{ \{(x,y) \mid x,y\in L, \, xy=c\} \mid c\in
L\}.
\end{array} \]

Then, $(\mathcal P, \mathcal L = \mathcal L_1 \cup \mathcal L_2
\cup \mathcal L_3 )$ is a $3$-net. The lines of these classes are
also called {\em horizontal, vertical} and {\em transversal
lines,} respectively. The point $O=(e,e)$ is the {\em origin} of
the net.

Let us now consider a $3$-net $(\mathcal P, \mathcal L = \mathcal
L_1 \cup \mathcal L_2 \cup \mathcal L_3 )$. Let $O\in \mathcal P$
be an arbitrary point, and let $\ell$, $k$ be the unique
horizontal and vertical lines through $O$, respectively. Then the
construction of Figure~1 defines a loop operation on $\ell$ with
neutral element $O$.

 \setlength{\unitlength}{2.05pt}
\begin{center}
\begin{picture}(100,42)  
\put(45,30){\line(5,-3){35}} \put(25,30){\line(5,-3){35}}
\put(25,30){\line(1,0){20}}  \put(45,9){\line(0,1){21}}
\thicklines \put(25,4){\line(0,1){35}}  
\put(20,9){\line(1,0){70}} \thinlines \put(45,30){\circle*{2}
$(x,y)$} \put(25,9){\circle*{2}}\put(44,5){$x$} \put(58,5){$y$}
\put(76,5){$x\cdot y$} \put(92,7.5){$\ell$}
 \put(21,29.5){$y$} \put(24.3,41.7){$k$} \put(1,3){$(e,e)=O$}
 \end{picture}

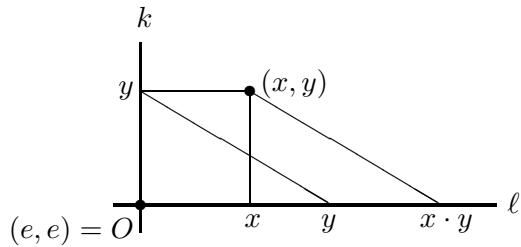
{Figure 1: The geometric definition of the coordinate loop.}
\end{center}

Since the parallel projections are bijections between lines of
different type, we can index the points of $k$ by points of
$\ell$, thus obtaining a bijection between $\mathcal P$ and $\ell
\times \ell$. The three line classes are determined by the
equations $X=c$, $Y=c$, $XY=c$, respectively, where $c$ is a
constant. We say that $(\ell,O)$ is a \emph{coordinate loop} of
the $3$-net $(\mathcal P$, $\mathcal L)$.

\subsection{Collineations and autotopisms}\label{Ss:Coli}

Let $\mathcal N = (\mathcal P, \mathcal L)$ be a $3$-net. {\em
Collineations} are line preserving bijective maps $\mathcal P \to
\mathcal P$. The group of collineations of $\mathcal N$ is denoted
by $\Coll\mathcal N$.  A collineation induces a permutation of the
line classes. There is therefore a group homomorphism from
$\Coll\mathcal N$ to the symmetric group $S_3$. The kernel of this
homomorphism consists of the {\em direction preserving
collineations.}

Let $L$ be the coordinate loop of $\mathcal N = (\mathcal
P,\,\mathcal L)$ with respect to some origin $O\in \mathcal P$.
Let $\varphi:\mathcal P \to \mathcal P$ be a bijection. Then
$\varphi$ preserves the line classes $1$ and $2$ if and only if it
has the form $(x,y)\mapsto (x\alpha, y\beta)$ for some bijections
$\alpha$, $\beta: L \to L$. Moreover, if $\varphi$ preserves the
line classes $1$ and $2$ then $\varphi$ also preserves the third
class if and only if there is a bijection $\gamma:L\to L$ such
that the triple $(\alpha,\beta,\gamma)$ is an autotopism of $L$.
Automorphisms of $L$ can be characterized in a similar way (see
Lemma \ref{Lm:GeomChar}).

\subsection{Bol reflections}

Let $\mathcal N$ be a $3$-net and $\ell_i\in\mathcal L_i$, for
some $i$. We define a certain permutation $\sigma_{\ell_i}$ on the
point set $\mathcal P$ (cf.\ Figure 2). For $P\in\mathcal P$, let
$a_j$ and $a_k$ be the lines through $P$ such that $a_j\in\mathcal
L_j$, $a_k\in\mathcal L_k$, and $\{i,\,j,\,k\}=\{1,\,2,\,3\}$.
Then there are unique intersection points $Q_j=a_j\cap \ell_i$,
$Q_k=a_k\cap \ell_i$. We define $P\sigma_{\ell_i}=b_j\cap b_k$,
where $b_j$ is the unique $j$-line through $Q_k$, and $b_k$ the
unique $k$-line through $Q_j$. The permutation $\sigma_{\ell_i}$
is clearly an involution satisfying $\mathcal
L_j\sigma_{\ell_i}=\mathcal L_k$, $\mathcal
L_k\sigma_{\ell_i}=\mathcal L_j$. If it happens to be the case
that $\sigma_{\ell_i}$ is a collineation, we call it the \emph{Bol
reflection with axis $\ell_i$}.

\setlength{\unitlength}{0.81mm}  
\begin{center}
\begin{picture}(120,37)(-3,10)
\thicklines \put(50,10){\line(0,1){35}} \thinlines
\put(25,10){\line(0,1){20}}
\multiput(75,30)(0,4){4}{\line(0,1){3}} \put(25,20){\circle*{1.7}}
\put(50,20){\circle*{1.7}}
\put(50,35){\circle*{1.7}} \put(75,35){\circle*{1.7}} 
\put(22,18){\line(5,3){31}} \put(78,37){\line(-5,-3){31}}
\put(22,20){\line(1,0){31}} \put(47,35){\line(1,0){31}}
\put(20,21){$P$} \put(79,37){$P^\prime = {P\sigma_{\ell_i}}$}
\put(51,10){$\ell_i$} \put(40,22){$a_k$} \put(35,29){$a_j$}
\put(60,37){$b_k$} \put(65,25){$b_j$} \put(53,17){$Q_k$}
\put(53,38){$Q_j$}
\end{picture}

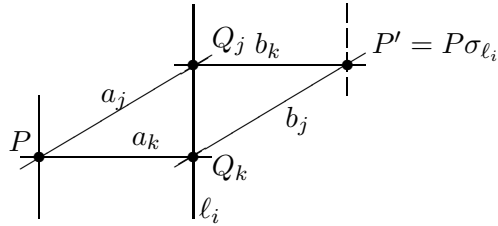
\vspace{6pt} {Figure 2: The Bol reflection with axis $\ell_i$.}
\end{center}

Obviously, every Bol reflection fixes a line pointwise (namely its
axis) and interchanges the other two line classes. In fact, it is
easy to see that any collineation with this property is a Bol
reflection. Then for any $\gamma\in\Coll{\mathcal N}$ and
$\ell\in\mathcal L$ we must have
$\gamma^{-1}\sigma_\ell\gamma=\sigma_{\ell\gamma}$, as
$\gamma^{-1}\sigma_{\ell}\gamma$ is a collineation fixing the line
$\ell\gamma$ pointwise. In words, the set of Bol reflections of
$\mathcal N$ is invariant under conjugations by elements of the
collineation group of $\mathcal N$.

Let $\ell_i \in \mathcal L_i$, $i=1,2,3$, be the lines through
some point $P$ of $\mathcal N$. As we have just seen,
$\sigma_{\ell_1} \sigma_{\ell_2} \sigma_{\ell_1} =
\sigma_{\ell_3}$, since $\ell_3\sigma_{\ell_1}=\ell_2$. Therefore
$(\sigma_{\ell_1} \sigma_{\ell_2})^3 = \id$ and
$\langle\sigma_{\ell_1}$, $\sigma_{\ell_2}$,
$\sigma_{\ell_3}\rangle$ is isomorphic to $S_3$. This fact will be
of importance later.

A $3$-net $\mathcal N$ is called a \emph{Moufang $3$-net} if
$\sigma_{\ell}$ is a Bol reflection for every line $\ell$. The
terminology is justified by Bol, who proved that $\mathcal N$ is a
Moufang $3$-net if and only if all coordinate loops of $\mathcal
N$ are Moufang \cite[p.\ 120]{Bruck}.

\setlength{\unitlength}{0.45mm}
\begin{center}
\begin{picture}(110,65)(-8,-7)
\thicklines \put(50,0){\line(0,1){60}} \thinlines
\put(0,0){\line(0,1){35}} \put(0,0){\line(1,0){50}}
\put(0,0){\line(2,1){50}} \put(0,35){\line(1,0){50}}
\put(0,35){\line(2,1){50}} \put(50,0){\line(2,1){50}}
\put(50,25){\line(1,0){50}} \put(50,35){\line(2,1){50}}
\put(50,60){\line(1,0){50}} \put(100,25){\line(0,1){5}}
\put(100,35){\line(0,1){5}} \put(100,45){\line(0,1){5}}
\put(100,55){\line(0,1){5}} \put(-8,-3){$P$} \put(-8,32){$R$}
\put(49,-7){$\ell$} \put(102,22){$P\sigma_\ell$}
\put(102,57){$R\sigma_\ell$} \put(0,0){\circle*{3}}
\put(50,0){\circle*{3}} \put(50,25){\circle*{3}}
\put(50,35){\circle*{3}} \put(0,35){\circle*{3}}
\put(50,60){\circle*{3}} \put(100,25){\circle*{3}}
\put(100,60){\circle*{3}}
\end{picture}

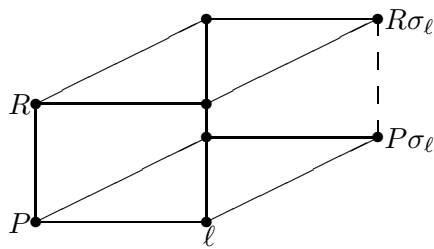
Figure 3: The $2$-Bol configuration.
\end{center}

The configuration in Figure 3 is called the \emph{$2$-Bol
configuration}. Using the other two directions of axes, we obtain
$1$- and $3$-Bol configurations. With these configurations at
hand, we see that the net $\mathcal N$ is Moufang if and only if
all its Bol configurations close (i.e., $R\sigma_\ell$ and
$P\sigma_\ell$ are collinear). See \cite[Sec.\ II.3]{Pflug} for
more on closures of net configurations.

\section{Composition algebras}

The most famous nonassociative Moufang loop is the multiplicative
loop of real octonions. Recall that octonions are built up from
quaternions in a way analogous to the construction of quaternions
from complex numbers, or complex numbers from real numbers.
Following Springer and Veldkamp \cite{SpringerVeldkamp}, we will
imitate this procedure over any field. We then construct a
countable family of finite simple Moufang loops, one for every
finite field $GF(q)$.

Let $F$ be a field and $V$ a vector space over $F$. A map $N:V\to
F$ is a \emph{quadratic form} if
$\langle\phantom{x},\phantom{x}\rangle : V\times V\to F$ defined
by $\langle u,v\rangle = (u+v)N-uN-vN$ is a bilinear form, and if
$(\lambda u)N =\lambda^2(uN)$ holds for every $u\in V$ and
$\lambda\in F$.

When $f:V\times V\to F$ is a bilinear form, then $u$, $v\in V$ are
\emph{orthogonal} (with respect to $f$) if $(u,v)f=0$. We write
$u\perp v$. The \emph{orthogonal complement} $W^\perp$ of a
subspace $W\le V$ is the subspace $\{v\in V;\;v\perp w$ for every
$w\in W\}$. The bilinear form $f$ is said to be
\emph{non-degenerate} if $V^\perp=0$. A quadratic form $N$ is
\emph{non-degenerate} if the bilinear form
$\langle\phantom{x},\phantom{x}\rangle$ associated with $N$ is
non-degenerate. When $N$ is non-degenerate, the vector space $V$
is said to be \emph{nonsingular}. A subspace $W$ of $(V,N)$ is
\emph{totally isotropic} if $uN=0$ for every $u\in W$. All maximal
totally isotropic subspaces of $(V,N)$ have the same dimension,
called the \emph{Witt index}. If $N$ is non-degenerate and $\dim
V\le \infty$ then the Witt index cannot exceed $\dim V/2$.

In this paper, an \emph{algebra} over $F$ is a vector space over
$F$ with bilinear multiplication. Specifically, we do not assume
that multiplication in an algebra is associative.

A \emph{composition algebra} $C=(C,N)$ over $F$ is an algebra with
a multiplicative neutral element $e$ such that the quadratic form
$N:C\to F$ is non-degenerate and
\begin{equation}\label{Eq:PreservesMult}
    (uv)N=uNvN
\end{equation}
holds for every $u$, $v\in C$. In this context, the quadratic form
$N$ is called a \emph{norm}.

When $\langle\phantom{x},\phantom{y}\rangle$ is the bilinear form
associated with the norm $N$, the \emph{conjugate} of $x\in C$ is
the element $\ol{x}=\langle x,e\rangle e - x$. Every element $x\in
C$ satisfies
\begin{displaymath}
    x^2-\langle x,e\rangle x +(xN)e=0
\end{displaymath}
(cf.\ \cite[Prop.\ 1.2.3]{SpringerVeldkamp}), and thus also
$x\overline{x}=\overline{x}x=(xN)e$. In particular, the
multiplicative inverse of $x$ is $x^{-1}=(xN)^{-1}\overline{x}$,
as long as $xN\ne 0$. Furthermore, $0\ne x\in C$ is a zero divisor
if and only if $xN=0$.

\subsection{The Cayley-Dickson process}

Let $C=(C,N)$ be a composition algebra over $F$ and $\lambda\in
F^*=F\setminus\{0\}$. Define a new product on $D=C\times C$ by
\begin{displaymath}
    (x,\,y)(u,\,v)=(xu+\lambda\ol{v}y,\,vx+y\ol{u}),
\end{displaymath}
where $x$, $y$, $u$, $v$ are elements of $C$. Also define the norm
$M$ on $D$ by
\begin{displaymath}
    (x,\,y)M = xN-\lambda(yN),
\end{displaymath}
where $x$, $y\in C$. By \cite[Prop.\ 1.5.3]{SpringerVeldkamp}, if
$C$ is associative then $D=(D,M)$ is a composition algebra.
Moreover, $D$ is associative if and only if $C$ is commutative and
associative. The above procedure is known as the
\emph{Cayley-Dickson process}.

We would now like to construct all composition algebras by
iterating the Cayley-Dickson process starting with $F$. However,
there is a twist when $F$ is of characteristic $2$. Namely, when
$\chr{F}=2$ then $F$ is not a composition algebra since $\langle
x,x\rangle = (x+x)N-xN-xN =0$ for every $x\in F$, thus $\langle
x,y\rangle = \langle x,\lambda x\rangle = \lambda\langle x,
x\rangle = 0$ for every $x$, $y\in F$, and $N$ is therefore
degenerate. The situation looks as follows:

\begin{theorem}[Thm.\ 1.6.2. \cite{SpringerVeldkamp}] Every
composition algebra over $F$ is obtained by iterating the
Cayley-Dickson process, starting from $F$ if $\chr{F}$ is not
equal to $2$, and from a $2$-dimensional composition algebra when
$\chr{F}$ is equal to $2$. The possible dimensions of a
composition algebra are $1$, $2$, $4$ and $8$. Composition
algebras of dimension $1$ or $2$ are commutative and associative,
those of dimension $4$ are associative but not commutative, and
those of dimension $8$ are neither commutative nor associative.

A composition algebra of dimension $2$ over $F$ is either a
quadratic field extension of $F$ or is isomorphic to $F\oplus F$.
\end{theorem}

For a generalization of composition algebras into dimension $16$
we refer the reader to \cite{JDH}.

\subsection{Split octonion algebras}

Composition algebras of dimension $8$ are known as \emph{octonion
algebras}. Since there is a parameter $\lambda$ in the
Cayley-Dickson process, it is conceivable (and sometimes true)
that there exist two octonion algebras over $F$ that are not
isomorphic.

A composition algebra $(C,N)$ is called \emph{split} if there is
$0\ne x\in C$ such that $xN=0$. By \cite[Thm.\
1.8.1]{SpringerVeldkamp}, over any field $F$ there is exactly one
split composition algebra in dimension $2$, $4$ and $8$, up to
isomorphism. As we have already noticed, split composition algebras
are precisely composition algebras with zero divisors. The unique
split octonion algebra over $F$ will be denoted by $\octo{F}$. (It
is worth mentioning that when $F$ is finite then every octonion
algebra over $F$ is isomorphic to $\octo{F}$, cf. \cite[p.\
22]{SpringerVeldkamp}.)

All split octonion algebras $\octo{F}$ were known already to Zorn,
who constructed them using the \emph{vector matrices}
\begin{equation}\label{Eq:Element}
    x=\vm{a}{\alpha}{\beta}{b},
\end{equation}
where $a$, $b\in F$ and $\alpha$, $\beta$ are vectors in $F^3$.
The norm $N$ is given as the ``determinant''
$\det{x}=ab-\alpha\cdot\beta$, where $\alpha\cdot\beta$ is the
usual dot product
\begin{displaymath}
    (\alpha_1,\alpha_2,\alpha_3)\cdot (\beta_1,\beta_2,\beta_3)
    = \alpha_1\beta_1+\alpha_2\beta_2+\alpha_3\beta_3.
\end{displaymath}
 The conjugate of $x$ is
\begin{equation}\label{Eq:Conjugate}
    \overline{x}=\vm{b}{-\alpha}{-\beta}{a},
\end{equation}
and two vector matrices are multiplied according to
\begin{equation}\label{Eq:Zorn}
    \vm{a}{\alpha}{\beta}{b}\vm{c}{\gamma}{\delta}{d} =
    \vm{ac+\alpha\cdot\delta}{a\gamma+d\alpha-\beta\times\delta}
    {c\beta+b\delta+\alpha\times\gamma}{\beta\cdot\gamma+bd},
\end{equation}
where $\beta\times\delta$ is the usual vector product
\begin{displaymath}
    (\beta_1,\beta_2,\beta_3)\times(\delta_1,\delta_2,\delta_3)
    =(\beta_2\delta_3-\beta_3\delta_2,\,
    \beta_3\delta_1-\beta_1\delta_3,\,
    \beta_1\delta_2-\beta_2\delta_1).
\end{displaymath}

The reader can think of this \emph{Zorn vector algebra} anytime we
speak of $\octo{F}$.

It turns out that every composition algebra satisfies the
\emph{alternative laws}
\begin{displaymath}
    (xy)x = x(yx),\quad x(xy)=x^2y, \quad (xy)y = xy^2.
\end{displaymath}
This is an easy corollary of the (not so easy) fact that
composition algebras satisfy the Moufang identities
$(\ref{Eq:Moufang})$, cf.\ \cite[Prop.\ 1.4.1]{SpringerVeldkamp}.

\section{A class of classical simple Moufang loops}

\subsection{Paige loops}

Although the octonion algebra $\octo{F}$ satisfies the Moufang
identities, it is not a Moufang loop yet, since it is not even a
quasigroup ($0\cdot x = 0$ for every $x\in\octo{F}$). Denote by
$M(F)$ the subset of $\octo{F}$ consisting of all elements of norm
(determinant) $1$. We have $\det{x}\det{y}=\det{xy}$ since
$\octo{F}$ is a composition algebra, which means that $M(F)$ is
closed under multiplication. The neutral element of $M(F)$ is
\begin{displaymath}
    e=\vm{1}{(0,0,0)}{(0,0,0)}{1},
\end{displaymath}
and the two-sided inverse of $x\in M(F)$ is $x^{-1}=\overline{x}$,
where $x$ is as in \eqref{Eq:Element} and $\overline{x}$ is as in
\eqref{Eq:Conjugate}.

Let $Z$ be the center of the Moufang loop $M(F)$. We have
$Z=\{e\}$ when $\chr{F}=2$, and $Z=\{e,-e\}$ when $\chr{F}\ne 2$.
Denote by $\paige{F}$ the Moufang loop $M(F)/Z$.

\begin{theorem}[Paige \cite{Paige}]\label{Th:Paige}
Let $F$ be a field and $M^*(F)$ the loop of Zorn vector matrices
of norm one modulo the center, multiplied according to
$(\ref{Eq:Zorn})$. Then $\paige{F}$ is a nonassociative simple
Moufang loop. When $F=GF(q)$ is finite, the order of $\paige{F}$
is $\frac{1}{d}q^3(q^4-1)$, where $d=(2,q-1)$.
\end{theorem}

The noncommutative loops $M^*(F)$ of Theorem \ref{Th:Paige} are
sometimes called \emph{Paige loops}.

In the remaining part of this section, we investigate the
multiplication groups of loops $M(F)$ and $M^*(F)$ constructed
over an arbitrary field $F$.

\subsection{Orthogonal groups}

Let $V$ be a vector space over $F$ with a non-degenerate quadratic
form $N:V\to F$. A linear transformation $f:V \to V$ is {\em
orthogonal with respect to $N$} if it preserves $N$, i.e., if
$(xf)N=xN$ for all $x \in V$. Then $f$ preserves the associated
bilinear form $\langle\phantom{x},\phantom{x}\rangle$ as well:
\begin{eqnarray*}
    \langle xf,yf \rangle &=& (xf+yf)N-(xf)N-(yf)N \\
    &=& (x+y)N-(x)N-(y)N \\
    &=& \langle x,y\rangle.
\end{eqnarray*}
The group consisting of all orthogonal transformations of $(V,N)$
is known as the {\em orthogonal group} $O(V)=O(V,N)$. The
determinant of an orthogonal transformation is $\pm 1$. Orthogonal
transformations with determinant $1$ form the {\em special
orthogonal group} $SO(V)$. The elements of $SO(V)$ are called
\emph{rotations}. One usually denotes by $\Omega(V)$ the
commutator subgroup $O'(V)$ of $O(V)$. By definition, every
element of $\Omega(V)$ is a rotation, and we would like to see
which rotations belong to $\Omega(V)$.

Take an element $g \in SO(V)$ and consider the map $1-g:x\mapsto
x-xg$. Define the bilinear form $\chi_g$ on $V(1-g)$ by
$(u,v)\chi_g=\langle u, w \rangle$, where $w$ is an arbitrary
vector from $V$ satisfying $w(1-g)=v$. Then $\chi_g$ is
well-defined and non-degenerate, by \cite[Thm.\ 11.32]{Taylor}.
Recall that the \emph{discriminant} $\discr(f)$ of a bilinear form
$f$ with respect to some basis is the determinant of its matrix.
Whether or not the discriminant of $\chi_g$ is a square in $F$
does not depend on the choice of the basis in $V(1-g)$. This
property characterizes elements of $\Omega(V)$.

\begin{lemma}[11.50 Thm.\ \cite{Taylor}] \label{lem:omega}
The rotation $g\in SO(V)$ belongs to $\Omega(V)$ if and only if
$\discr(\chi_g) \in F^2$.
\end{lemma}

Pick any element $\sigma \in O(V)$ with $\sigma^2=\id$. The two
subspaces
\begin{eqnarray*}
    U&=&V(\sigma -1) = \{ x\sigma -x \mid x \in V \},\\
    W&=&V(\sigma +1) = \{x\sigma +x \mid x \in V \}
\end{eqnarray*}
are orthogonal to each other. Indeed,
\begin{displaymath}
    \langle x\sigma -x, y \sigma + y \rangle = \langle x\sigma , y \rangle
    - \langle x, y \sigma \rangle = \langle x\sigma , y \rangle - \langle
    x\sigma , y \sigma^2 \rangle = 0.
\end{displaymath}

The subspace $W$ consists of vectors invariant under $\sigma$. If
$W$ is a nonsingular hyperplane (that is, a subspace of dimension
$\dim V -1$) then $\sigma$ is called a {\em symmetry with respect
to $W$.} (If $\chr(F)=2$ then $\sigma$ is usually called a {\em
transvection}.) If $\sigma$ is a symmetry with respect to $W$ and
$g\in O(V)$, the conjugate $\sigma^g=g^{-1}\sigma g$ is a symmetry
with respect to $Wg$.

\subsection{Multiplication groups of Paige loops}

Let now $V=\octo{F}$ be the split octonion algebra over $F$. We
identify the vector matrix
\[ x = \vm{x_0}{(x_1,x_2,x_3)}{(x_4,x_5,x_6)}{x_7} \]
with the column vector $(x_0,\ldots,x_7)^t$, and we use the
canonical basis of $F^8$ as the basis of $V$. Since $\langle
x,y\rangle = \det(x+y)-\det x-\det y = x_7y_0-x_4y_1 - x_5y_2 -
x_6y_3 - x_1y_4 - x_2y_5 - x_3y_6 + x_0y_7$, the bilinear from
$\langle x,y\rangle$ can be expressed as $x^t J y$, where
\begin{equation}\label{Eq:Matrix}
J = \left ( \begin{array}{cccccccc}
0  & 0 & 0 & 0 & 0 & 0 & 0 &  1 \\
0  & 0 & 0 & 0 & -1 & 0 & 0 &  0 \\
0  & 0 & 0 & 0 & 0 & -1 & 0 &  0 \\
0  & 0 & 0 & 0 & 0 & 0 & -1 &  0 \\
0  & -1 & 0 & 0 & 0 & 0 & 0 &  0 \\
0  & 0 & -1 & 0 & 0 & 0 & 0 &  0 \\
0  & 0 & 0 & -1 & 0 & 0 & 0 &  0 \\
1  & 0 & 0 & 0 & 0 & 0 & 0 &  0
       \end{array} \right ).
\end{equation}

Recall that $M(F)$ consists of all elements of $\octo{F}$ that are of norm $1$.

\begin{lemma} \label{lem:mult1}
For every $a \in M(F)$, we have $L_a,R_a \in \Omega(V)$.
\end{lemma}
\begin{proof}
We only deal with the case $L_a$. Since $aN=1$, we have $L_a\in
O(V)$, by \eqref{Eq:PreservesMult}. Let $a=(a_0$, $\dots$,
$a_7)^t$ and write matrix maps to the left of their arguments.
Then $L_a$ can be identified with
\[ \left ( \begin{array}{cccccccc}
 a_0&   0&   0&   0& a_1& a_2& a_3&   0 \\
   0& a_0&   0&   0&   0& a_6&-a_5& a_1 \\
   0&   0& a_0&   0&-a_6&   0& a_4& a_2 \\
   0&   0&   0& a_0& a_5&-a_4&   0& a_3 \\
 a_4&   0&-a_3& a_2& a_7&   0&   0&   0 \\
 a_5& a_3&   0&-a_1&   0& a_7&   0&   0 \\
 a_6&-a_2& a_1&   0&   0&   0& a_7&   0 \\
   0& a_4& a_5& a_6&   0&   0&   0& a_7
       \end{array} \right ). \]

Routine calculation yields $\det(L_a)=(aN)^4$, and $L_a \in SO(V)$
follows. By Lemma \ref{lem:omega}, it suffices to show
$\discr(\chi_{L_a}) \in F^2$.

Assume first that $(e-a)N\ne 0$. Then $V(1-L_a)=V(e-a)=V$, and
$((e-a)^{-1}v)(1-L_a)=v$ for every $v\in V$. Thus $(u,v)\chi_{La}
= \langle u,vL_{e-a}^{-1}\rangle$, and the matrix of $\chi_{L_a}$
is $JL_{e-a}^{-1}$, where $J$ is as in \eqref{Eq:Matrix}.
Therefore $\discr(\chi_{L_a}) = \det(J)\det(L_{e-a})^{-1} =
((e-a)N)^{-4}\in F^2$.

Suppose now $(e-a)N=0$ and exclude the trivial case $e-a \in F$.
Define the elements
\[ \begin{array}{l}
e_0 = \vm{1}{(0,0,0)}{(0,0,0)}{0},\quad
e_1 = \vm{0}{(1,0,0)}{(0,0,0)}{0}, \\[12pt]
e_2 = \vm{0}{(0,1,0)}{(0,0,0)}{0},\quad e_3 =
\vm{0}{(0,0,1)}{(0,0,0)}{0}
   \end{array} \]
and
\[ \begin{array}{l}
f_0= (e-a)e_0 = \vm{1-a_0}{(0,0,0)}{(-a_4,-a_5,-a_6)}{0}, \\[12pt]
f_1= (e-a)e_1 = \vm{0}{(1-a_0,0,0)}{(0,-a_3,a_2)}{-a_4}, \\[12pt]
f_2= (e-a)e_2 = \vm{0}{(0,1-a_0,0)}{(a_3,0,-a_1)}{-a_5}, \\[12pt]
f_3= (e-a)e_3 = \vm{0}{(0,0,1-a_0)}{(-a_2,a_1,0)}{-a_6}.
   \end{array} \]
The vectors $e_i$ span a totally isotropic subspace of $V$ and
$f_i \in (e-a)V$. Since $\langle (e-a)x,(e-a)y \rangle =
(e-a)N\langle x,y \rangle =0$, $(e-a)V$ is totally isotropic as
well. In particular, $\dim ((e-a)V) \leq 4$.

Assume $a_0 \neq 1$. Then, the vectors $f_i$ are linearly
independent and hence form a basis of $(e-a)V$. The matrix
$M=(m_{ij})$ of $\chi_{L_a}$ with respect to the basis $\{f_0$,
$f_1$, $f_2$, $f_3\}$ satisfies $m_{ij} = (f_i,f_j)\chi_{L_a} =
\langle f_i, e_j \rangle$, which yields
\[ M= \left ( \begin{array}{cccc}
0 & a_4 & a_5 & a_6 \\
-a_4 & 0 & a_3 & -a_2 \\
-a_5 & -a_3 & 0 & a_1 \\
-a_6 & a_2 & -a_1 & 0 \\
          \end{array} \right ), \]
by calculation. Then $\discr(\chi_{L_a}) = \det{M} =
(a_1a_4+a_2a_5+a_3a_6)^2 \in F^2$.

The special case $a_0=1$ can be calculated similarly.
\end{proof}

For the rest of this section, let $\iota$ denote the conjugation
map $x\mapsto \overline{x}$. Note that $\iota\in O(V)$ and $e\iota
= e$.

\begin{lemma}
Any element $g\in O(V)$ with $eg=e$ commutes with $\iota$.
\end{lemma}
\begin{proof}
We have $\overline{x} g = (\langle x,e \rangle e - x)g = \langle
x,e \rangle eg - xg = \langle xg,eg \rangle eg - xg = \langle xg,e
\rangle e - xg = \overline{xg}$.
\end{proof}

\begin{lemma}
For an arbitrary element $g\in O(V)$, we define $\iota^g =
g^{-1}\iota g$. Put $a=eg$. Then $aN=1$ and $x\iota^g =
a\overline{x}a$ holds for all $x\in V$.
\end{lemma}
\begin{proof}
On the one hand, $aN=(eg)N=eN=1$, therefore $a\overline{a}=e$ and
$a^{-1}=\overline{a}$. On the other hand, $g=h L_a$ for some $h$
with $eh=e$. By the previous lemma, $\iota^g=L_a^{-1} \iota L_a$
and $x\iota^g = ((xL_a^{-1})\iota)L_a = a(\overline{a^{-1}x}) =
a\overline{x}a$.
\end{proof}

The map $-\iota:x \mapsto -\bar x$ is a symmetry with respect to
the $7$-dimensional nonsingular hyperplane
\[ H = \left \{ \vm{x_0}{(x_1,x_2,x_3)}{(x_4,x_5,x_6)}{-x_0} \Bigg|\ x_i \in F
\right \}. \] The conjugate $-\iota^g$ is a symmetry with respect
to $Hg$. This means that
\[\mathcal C = \{ -\iota^g \mid g \in O(V) \} = \{ -L_a^{-1} \iota L_a \mid a
\in M(F) \} \] is a complete conjugacy class consisting of
symmetries.

\begin{theorem} \label{thm:mult(mf)}
The multiplication group of $M(F)$ is $\Omega(\octo{F},N)$.
\end{theorem}
\begin{proof}
By Lemma \ref{lem:mult1}, $\Mlt(M(F)) \leq \Omega(V)$. We have
$(ax)\iota = x\iota \, \bar a$, which implies
\[\iota \iota^g = \iota L_a^{-1} \iota L_a = R_aL_a\]
and
\[\iota^g \iota^h = (\iota \iota^g)^{-1} (\iota \iota^h) = (R_aL_a)^{-1}
(R_bL_b) \in \Mlt(M(F))\]
 for $g,h \in O(V)$. Let us denote by
$\mathcal D$ the set consisting of $\iota^g \iota^h$, $g,h \in
O(V)$. $\mathcal D$ is clearly an invariant subset of $O(V)$. By
[1, Thm.\ 5.27], $\mathcal D$ generates $\Omega(V)$, which proves
$\Mlt(M(F)) = \Omega(\octo{F}, N)$.
\end{proof}

Finally, we determine the multiplication groups of Paige loops.

\begin{corollary}
The multiplication group of the Paige loop $M^*(F)$ is the simple
group $P\Omega(\octo{F},N) = P\Omega^+_8(F)$.
\end{corollary}
\begin{proof}
The surjective homomorphism $\varphi: M(F) \to M^*(F)$, $x \mapsto
\pm x$ induces a surjective homomorphism $\Phi:\Mlt(M(F)) \to
\Mlt(M^*(F))$. On the one hand, the kernel of $\Phi$ contains $\pm
\id$. On the other hand $P\Omega(V) = \Omega(V)/\{\pm\id\}$ is a
simple group, cf. [1, Thm. 5.27]. Since $\Mlt(M^*(F))$ is not
trivial, we must have $\Mlt(M^*(F)) = P\Omega(V)$. Finally, the
norm $N$ has maximal Witt index $4$, therefore the notation
$P\Omega^+_8(F)$ is justified.
\end{proof}

\begin{remark}{\rm
The result of Theorem \ref{thm:mult(mf)} is {\em folklore,} that
is, most of the authors (Freudenthal, Doro, Liebeck, etc.) use it
as a {\em well-known} fact without making a reference. The authors
of the present paper are not aware of any reference, however,
especially one that would handle all fields at once.}
\end{remark}

\section{Groups with triality}

\subsection{Triality}

Let $G$ be a group. We use the usual notation $x^y=y^{-1}xy$ and
$[x,y]=x^{-1}y^{-1}xy=x^{-1}x^y$ for $x$, $y\in G$. Let $\alpha$
be an automorphism of $G$, then $x\alpha$ will be denoted by
$x^\alpha$ as well, and $[x,\alpha]$ will stand for
$x^{-1}x^\alpha$. The element $\alpha^y \in \Aut{G}$ maps $x$ to
$x^{y^{-1}\alpha y} = ((x^{y^{-1}})^\alpha)^y$.

Let $S_n$ be the symmetric group on $\{1, \ldots, n\}$. $G \times
H$ and $G \rtimes H$ will stand for the direct and semidirect
product of $G$ and $H$, respectively. In the latter case, $H$ acts
on $G$.

We have the following definition due to Doro \cite{Doro}.

\begin{definition}
The pair $(G,S)$ is called a {\em group with triality}, if $G$ is
a group, $S \leq \Aut{G}$, $S=\langle \sigma,\rho|\;
\sigma^2=\rho^3=(\sigma\rho)^2=1 \rangle \cong S_3$, and for all
$g\in G$ the {\em triality identity}
\[ [g,\sigma]\,[g, \sigma]^\rho\, [g,\sigma]^{\rho^2}=1\]
holds.
\end{definition}

The principle of triality was introduced by Cartan \cite{Cartan}
in 1938 as a property of orthogonal groups in dimension $8$, and
his examples motivated Tits \cite{Tits}.  Doro was the first one
to define the concept of an abstract group with triality, away
from any context of a given geometric or algebraic object.

\begin{definition}
Let $(G_i,\langle \sigma_i, \rho_i \rangle)$, $i=1,2$ be groups
with triality. The homomorphism $\varphi:G_1 \to G_2$ is an {\em
$S$-homomorphism} if $g\sigma_1\varphi = g\varphi\sigma_2$ and
$g\rho_1\varphi = g\varphi\rho_2$ hold for all $g\in G_1$. The
kernel of an $S$-homomorphism is an $S$-invariant normal subgroup.
The group with triality $G$ is said to be {\em $S$-simple} if it
has no proper $S$-invariant normal subgroups.
\end{definition}

The following examples of groups with triality are of fundamental
importance. They are adopted from Doro \cite{Doro}.

\begin{example} \label{ex:typeI}
Let $A$ be a group, $G=A^3$, and let $\sigma$, $\rho\in\Aut{G}$ be
defined by $\sigma:(a_1,a_2,a_3) \mapsto (a_2,a_1,a_3)$ and $\rho:
(a_1,a_2,a_3) \mapsto (a_2,a_3,a_1)$. Then $G$ is a group with
triality with respect to $S=\langle \sigma, \rho \rangle$.
\end{example}
\begin{example} \label{ex:typeII}
Let $A$ be a group with $\varphi \in \Aut{A}$,
$\varphi\ne\mathrm{id}_A$, satisfying $x \, x^{\varphi} \,
x^{\varphi^2}=1$ for all $x \in A$. Put $G=A \times A$, $\sigma:
(a_1,a_2) \mapsto (a_2,a_1)$ and $\rho: (a_1,a_2) \mapsto
(a_1^\varphi, a_2^{\varphi^{-1}})$. Then $G$ is a group with
triality with respect to $S=\langle \sigma, \rho \rangle$.

If $A$ is of exponent $3$ and $\varphi=\mathrm{id}_A$ then $G$ is
a group with triality in a wider sense, meaning that the triality
identity is satisfied but $S$ is not isomorphic to $S_3$.
\end{example}
\begin{example} \label{ex:comm}
Let $V$ be a two-dimensional vector space over a field of
characteristic different from $3$. Let $S$ be the linear group
generated by the  matrices
\[ \rho = \vm{-1}{-1}{1}{0} \hskip 1cm \mbox{and} \hskip 1cm
\sigma= \vm{0}{1}{1}{0}. \] Then the additive group of $V$ and $S$
form a group with triality.
\end{example}
\begin{remark}\rm
If $A$ is a simple group then the constructions in Examples
\ref{ex:typeI} and \ref{ex:typeII} yield $S$-simple groups with
triality. Obviously, if $(G,S)$ is a group with triality and $G$
is simple (as a group) then $(G,S)$ is $S$-simple. Below, we are
concerned with the converse of this statement.
\end{remark}

\subsection{Triality of Moufang nets}

In the following, $(G,S)$ stands for a group $G$ with automorphism
group $S$ isomorphic to $S_3$. Let $\sigma$, $\rho\in S$ be such
that $\sigma^2=\rho^3=\mathrm{id}$. Let the three involutions of
$S$ be $\sigma_1=\sigma$, $\sigma_2=\sigma \rho$ and
$\sigma_3=\rho\sigma =\sigma \rho^2 $. Finally, the conjugacy
class $\sigma_i^G$ will be denoted by $\mathcal C_i$.

The following lemma gives a more conceptual reformulation of
Doro's triality. (It is similar to Lemma 3.2 of \cite{Liebeck},
attributed by Liebeck to Richard Parker.)

\begin{lemma} \label{l322}
The pair $(G,S)$ is a group with triality if and only if
$(\tau_i\tau_j)^3=\mathrm{id}$ for every $\tau_i\in\mathcal C_i$,
$\tau_j\in\mathcal C_j$, where $i$, $j\in\{1$, $2$, $3\}$ and
$i\ne j$. In this case, $(G,\langle\tau_i,\tau_j\rangle)$ is a
group with triality as well.
\end{lemma}
\begin{proof}
The condition of the first statement claims something about the
conjugacy classes ${\mathcal C}_i$, which do not change if we
replace $S$ by $\langle \tau_i, \tau_j \rangle$. This means that
the first statement implies the second one.

For the first statement, it suffices to investigate the case
$i=1$, $j=3$, $\tau_1=\sigma^g$ and $\tau_3=\sigma \rho^2$, with
arbitrary $g \in G$. Then the following equations are equivalent
for every $g\in G$:
\begin{eqnarray*}
    1&=&(\sigma^g (\sigma\rho^2))^3,\\
    1&=&\sigma^g (\sigma\rho^2) \cdot \sigma^g(\sigma\rho^2) \cdot
    \sigma^g (\sigma\rho^2),\\
    1&=& [g,\sigma]\rho^2 \cdot [g,\sigma]\rho^2 \cdot
    [g,\sigma]\rho^2,\\
    1&=&[g,\sigma]\cdot \rho^{-1}[g,\sigma]\rho \cdot \rho
    [g,\sigma]\rho^2,\\
    1&=&[g, \sigma]\,[g, \sigma]^{\rho}\, [g,\sigma]^{\rho^2}.
\end{eqnarray*}
This finishes the proof.
\end{proof}

The next lemma already foreshadows the relation between Moufang
$3$-nets and groups with triality.

\begin{lemma} \label{l323}
Let $P$ be a point of the Moufang $3$-net $\mathcal N$. Denote by
$\ell_1, \ell_2$ and $\ell_3$ the three lines through $P$ with
corresponding Bol reflections $\sigma_1$, $\sigma_2$, $\sigma_3$.
Then the collineation group $S=\langle \sigma_1, \sigma_2,
\sigma_3 \rangle\cong S_3$ acts faithfully on the set $\{\ell_1,
\ell_2, \ell_3\}$. This action is equivalent to the induced action
of $S$ on the parallel classes of $\mathcal N$.
\end{lemma}
\begin{proof}
As we have demonstrated in Section \ref{Sc:LN}, the conjugate of a
Bol reflection is a Bol reflection. Thus $\sigma_1 \sigma_2
\sigma_1=\sigma_3=\sigma_2 \sigma_1 \sigma_2$, which proves the
first statement. The rest is trivial.
\end{proof}

Using these lemmas, we can prove two key propositions.

\begin{proposition} \label{p324}
Let $\mathcal N$ be a Moufang $3$-net and let $M$ be the group of
collineations generated by all Bol reflections of $\mathcal N$.
Let $M_0\leq M$ be the direction preserving subgroup of $M$. Let
us fix an arbitrary point $P$ of $\mathcal N$ and denote by $S$
the group generated by the Bol reflections with axes through $P$.
Then $M_0\triangleleft M$, $M=M_0 S$, and the pair $(M_0,S)$ is a
group with triality.
\end{proposition}
\begin{proof}
$M_0\triangleleft M=M_0 S$ is obvious. Thus $S$ is a group of
automorphism of $M_0$ by conjugation. By Lemma \ref{l322}, it is
sufficient to show $\langle\sigma_i^g, \sigma_j^h \rangle \cong
S_3$ for all $g$, $h \in M_0$, where $\sigma_i$ and $\sigma_j$ are
the reflections on two different lines through $P$.  Since $g$,
$h$ preserve directions, the axes of $\sigma_i^g$ and $\sigma_j^h$
intersect in some point $P'$. Hence $\langle\sigma_1^g, \sigma_2^h
\rangle \cong S_3$, by Lemma \ref{l323}.
\end{proof}

The converse of the proposition is true as well.

\begin{proposition} \label{p325}
Let $(G,S)$ be a group with triality. The following construction
determines a Moufang $3$-net ${\mathcal N}(G,S)$. Let the three
line classes be the conjugacy classes ${\mathcal C}_1$, ${\mathcal
C}_2$ and ${\mathcal C}_3$. By definition, three mutually
non-parallel lines $\tau_i \in {\mathcal C}_i$ $(i=1,2,3)$
intersect in a common point if and only if
\[\langle \tau_1, \tau_2, \tau_3 \rangle \cong S_3.\]
Moreover, if $G_1=[G,S]S=\langle {\mathcal C}_1, {\mathcal C}_2,
{\mathcal C}_3 \rangle$, then the group $M({\mathcal N})$
generated by the Bol reflections of $\mathcal N$ is isomorphic to
$G_1/Z(G_1)$.
\end{proposition}
\begin{proof}
By definition, parallel lines do not intersect. When formulating
the triality identity as in Lemma \ref{l322}, we see that two
non-parallel lines have a point in common such that there is
precisely one line from the third parallel class incident with
this point. This shows that ${\mathcal N}(G,S)$ is a $3$-net
indeed.

The Moufang property follows from the construction immediately,
since one can naturally associate an involutorial collineation to
any line $\tau_i \in \mathcal C_i$, namely the one induced by
$\tau_i$ on $G$. This induced map $\bar\tau_i$ interchanges the
two other parallel classes $\mathcal C_j$, $\mathcal C_k$ and
fixes the points on its axis, that is, it normalizes the $S_3$
subgroups containing $\tau_i$.

Finally, since a Bol reflection acts on the line set in the same
way that the associated ${\mathcal C}_i$-element acts on the set
$\cup {\mathcal C}_j$ by conjugation, we have the isomorphism
$M({\mathcal N})\cong G_1/Z(G_1)$.
\end{proof}

\begin{remark}\rm
From the point of view of dual $3$-nets, the point set is the
union of the three classes ${\mathcal C}_i$, and lines consist of
the intersections of an $S_3$ subgroup with each of the three
classes.
\end{remark}
\begin{remark}\rm
One finds another construction of groups with triality using the
geometry of the associated $3$-net in P.~O.~Mikheev's paper
\cite{Mikheev}. A different approach to groups with triality is
given in J.~D.~Phillips' paper \cite{JD}.
\end{remark}

\subsection{Triality collineations in coordinates} \label{subsec:bolreflcoord}

At this point, we find it useful to write down the above maps in
the coordinate system of the $3$-net. If we denote by
$\sigma_m^{(v)}, \sigma_m^{(h)}, \sigma_m^{(t)}$ the Bol
reflections with axes $X=m$, $Y=m$, $XY=m$, respectively, then we
have
\begin{eqnarray*}
\sigma_m^{(v)}&:& (x,y) \mapsto (m (x^{-1}m), m^{-1}(xy)), \\[2pt]
\sigma_m^{(h)}&:& (x,y) \mapsto ((xy)m^{-1}, (my^{-1}) m), \\[2pt]
\sigma_m^{(t)}&:& (x,y) \mapsto (my^{-1},x^{-1} m).
\end{eqnarray*}
This yields
\begin{eqnarray*}
\sigma_m^{(v)} \sigma_1^{(v)} &:& (x,y) \mapsto (m^{-1} (xm^{-1}), my), \\[2pt]
\sigma_m^{(h)} \sigma_1^{(h)} &:& (x,y) \mapsto (xm, (m^{-1}y) m^{-1}), \\[2pt]
\sigma_{m^{-1}}^{(t)} \sigma_1^{(t)} &:& (x,y) \mapsto (mx, ym).
\end{eqnarray*}
These are direction preserving collineations generating $G$.

They can be written in the form $\sigma_m^{(v)} \sigma_1^{(v)} =
(L_m^{-1} R_m^{-1}, L_m)$, $\sigma_m^{(h)} \sigma_1^{(h)} = (R_m,
L_m^{-1} R_m^{-1})$ and $\sigma_m^{(t)} \sigma_1^{(t)} = (L_m,
R_m)$ as well. The associated autotopisms are
\begin{displaymath}
    (L_m^{-1} R_m^{-1}, L_m, L_m^{-1}),\quad
    (R_m, L_m^{-1}R_m^{-1}, R_m^{-1}),\quad
    (L_m, R_m, L_mR_m),
\end{displaymath}
respectively. By the way, the fact that these triples are
autotopisms is equivalent with the Moufang identities
\eqref{Eq:Moufang}.

\section{The classification of nonassociative finite simple Moufang loops}

\subsection{Simple $3$-nets}

The classification of finite simple Moufang loops is based on the
classification of finite simple groups with triality. Using the
results of the previous section, the classification can be done in
the following steps.

\begin{proposition}
Let $\varphi:{\mathcal N}_1\to {\mathcal N}_2$ be a map between
two $3$-nets that preserves incidence and directions.
\begin{enumerate}
\item[$(i)$] Suppose that $\varphi(P_1)=P_2$ holds for the points $P_1\in\mathcal
N_1$, $P_2\in\mathcal N_2$. Then $\varphi$ defines a homomorphism
$\bar\varphi:L_1\to L_2$ in a natural way, where $L_i$ is the
coordinate loop of the $3$-net ${\mathcal N}_i$ with origin $P_i$.
Conversely, the loop homomorphism $\bar\varphi:L_1\to L_2$
uniquely determines a collineation $\mathcal N_1\to\mathcal N_2$,
namely $\varphi$.

\item[$(ii)$] Suppose that the $3$-nets ${\mathcal N}_i$ $(i=1,2)$ are Moufang
and $\varphi$ is a collineation onto. Let us denote by $(M_i,S)$
the group with triality that corresponds to the $3$-net ${\mathcal
N}_i$. Then the maps $\sigma_\ell\mapsto \sigma_{\ell\varphi}$
induce a surjective $S$-homomorphism $\tilde\varphi:M_1\to M_2$,
where $\sigma_\ell$ is the Bol reflection in $\mathcal N_1$ with
axis $\ell$. Conversely, an $S$-homomorphism $M_1\to M_2$ defines
a direction preserving collineation between the $3$-nets
${\mathcal N}(M_1,S)$ and ${\mathcal N}(M_2,S)$.
\end{enumerate}
\end{proposition}
\begin{proof}
The first part of statement (i) follows from the geometric
definition of the loop operation in a coordinate loop; the second
part is trivial. For the (ii) statement, it is sufficient to see
that a relation of the reflections $\sigma_\ell$ corresponds to a
point-line configuration of the $3$-net, and that the
$\varphi$-image of the configuration induces the same relation on
the reflections $\sigma_{\varphi(\ell)}$. The converse follows
from Proposition \ref{p325}.
\end{proof}

In the sense of the proposition above, we can speak of {\em simple
$3$-nets}, that is, of $3$-nets having only trivial homomorphisms.
The next proposition follows immediately.

\begin{proposition} \label{all:egyszmoufh}
If $L$ is a simple Moufang loop, then the associated $3$-net
${\mathcal N}$ is simple as well. That is, the group $(M_0, S)$
with triality determined by ${\mathcal N}$ is $S$-simple.
\end{proposition}

\subsection{$S$-simple groups with triality}

The structure of $S$-simple groups with triality is rather
transparent.  It is clear that $G$ is $S$-simple if and only if it
has no $S$-invariant proper nontrivial normal subgroups.

Let $G$ be such a group and let $N\triangleleft G$ be an arbitrary
proper normal subgroup of $G$. Let us denote by $N_i$ the images
of $N$ under the elements of $S$, $i=1,\ldots,6$.

Since the union and the intersection of the groups $N_i$ is an
$S$-invariant normal subgroup of $G$, we have $G = N_1 \cdots N_6$
and $\{1\} = N_1 \cap \ldots \cap N_6$. If $N_i \cap N_j$ is a
proper subgroup of $N_i$ for some $i$, $j=1$, $\dots$, $6$, then
we replace $N_i$ by $N_i \cap N_j$. We can therefore assume that
the groups $N_i$ intersect pairwise trivially. Since $S$ acts
transitively on the groups $N_i$, one of the following cases must
occur:

\noindent{\bf Case A.} $G$ is a simple group. In this case, there
is no proper normal subgroup $N$.

\noindent{\bf Case B.} The number of distinct groups $N_i$ is $2$.
Then $N=N^\rho$, $M=N^\sigma$, $G=NM$, $N\cap M=\{1\}$ and
elements of $N$ and $M$ commute. Every element $g\in G$ can be
written as $g=ab^\sigma$. $\rho$ induces an automorphism $\varphi$
on $N$. Then, $g^\sigma = a^\sigma b = b a^\sigma$ and $g^\rho =
a^\rho b^{\sigma\rho} = a^\rho b^{\rho^{-1}\sigma} = a^\varphi
b^{\varphi^{-1}\sigma}$.

Moreover, applying the triality identity on $a \in N$, we obtain
\begin{displaymath}
    (a^{\varphi^2} a^\varphi a)^{-1} (aa^\varphi a^{\varphi^2})^\sigma = 1,
\end{displaymath}
which is equivalent with the identity $aa^\varphi a^{\varphi^2} =
1$. This means that the map $N\times N \to G$, $(a,b) \mapsto
ab^\sigma$ defines an $S$-isomorphism between $G$ and the
construction in Example \ref{ex:typeII}.

However, a result of Khukhro claims that the existence of the
automorphism $\varphi$ of $N$ implies that $N$ is nilpotent of
class at most $3$ (see \cite[p.\ 223]{Khukhro}, \cite[Thm.\
3.3]{NagyValsecchi}). Therefore, no $S$-simple group with triality
can be constructed in this case.

\noindent{\bf Case C.} The number of distinct groups $N_i$ is $3$:
$N=N_1$, $N^\rho=N_2$, $N^{\rho^2} = N_3$. We can assume
$N_1^\sigma=N_1$ and $N_2^\sigma=N_3$.

\noindent{\bf Case C/1.} Assume that $M=N_1\cap (N_2N_3)$ is a
proper subgroup of $N_1$. Then $M^\rho \in N_1^\rho = N_2
\subseteq N_2N_3$, similarly $M^{\rho^2} \in N_2N_3$. Moreover,
$M^\sigma = N_1^\sigma \cap (N_2^\sigma N_2^\sigma) = M$. Hence,
$MM^\rho M^{\rho^2}$ is a proper $S$-invariant normal subgroup of
$G$, a contradiction.

\noindent{\bf Case C/2.} Assume $N_1\cap (N_2N_3) = \{1\}$. Then
$G=N_1 \times N_2 \times N_3 \cong N^3$. By the triality identity,
we have $a^{-1}a^\sigma \in N_1 \cap (N_2N_3)$ for any $a\in N_1$,
thus, $a^\sigma=a$. Consider the map $\Phi:N^3 \to G$,
$\Phi(a,b,c)=ab^\rho c^{\rho^2}$. By
\[(ab^\rho c^{\rho^2})^\sigma = ac^\rho b^{\rho^2}
\hskip 1cm \mbox{and} \hskip 1cm (ab^\rho c^{\rho^2})^\rho =
ca^\rho b^{\rho^2},\] $\Phi$ defines an $S$-isomorphism between
$G$ and the group with triality in Example \ref{ex:typeI}.

\noindent{\bf Case C/3.} Assume $N_1 \subseteq N_2N_3$, $G$
noncommutative. We have $G=N_1 \times N_2 \cong N^2$. Since $G$ is
$S$-simple, we must have $Z(G) = \{1\}$ and $Z(N)=\{1\}$. Let us
assume that $a^\rho =a_1a_2$ with $1 \neq a_1 \in N_1$, $a_2 \in
N_2$ for some element $a \in N_1=N$. Take $b \in N$ with $a_1 b
\neq b a_1$. Every element of $N_1$ commutes with every element of
$N_2$. This implies
\[1 \neq [a_1a_2, b] = [a^\rho,b] \in N \cap N^\rho,\]
a contradiction.

\noindent{\bf Case C/4.} If $G$ is commutative and $S$-simple,
then we are in the situation of Example \ref{ex:comm}. The proof
is left to the reader.

We summarize these results in the following proposition. In the
finite case the result was proved by S.~Doro \cite{Doro}. In the
infinite case, it is due to G.~P.~Nagy and M.~Valsecchi
\cite{NagyValsecchi}.

\begin{proposition}
Let $G$ be a noncommutative $S$-simple group with triality. Then
either $G$ is simple or $G=A \times A \times A$, where $A$ is a
simple group and the triality automorphisms satisfy
$(a,b,c)^\rho=(c,a,b)$, $(a,b,c)^\sigma = (a,c,b)$.
\end{proposition}

\subsection{The classification}

\begin{lemma}
Let $G=A\times A \times A$ be an $S$-simple group with triality.
Then the associated loop is isomorphic to the group $A$.
\end{lemma}
\begin{proof}
We leave to the reader to check that an associative simple Moufang
loop $A$ has $G=A\times A \times A$ as a group with triality.
Since the group with triality determines the $3$-net uniquely, and
since groups are $G$-loops, that is, the coordinate loop does not
depend on the choice of the origin, we are done.
\end{proof}

Also the following result is due to Doro, but the way of proving
is based on the geometric approach, hence new.

\begin{lemma}[Doro]
Assume that $G$ is a group with triality such that $\rho$ is an
inner automorphism. Then the associated Moufang loop 
has exponent $3$.
\end{lemma}
\begin{proof}
Let $\rho$ be an inner automorphism of $G$. We assume $G$ to be a
group of direction preserving collineations of the $3$-net
$\mathcal N$ associated with $L$. We consider $\rho$ as a
collineation of $\mathcal N$ permuting the directions cyclicly. We
denote by $\Gamma^+$ the collineation group of $\mathcal N$
generated by $G$ and $\rho$. $\Gamma^+$ consists of collineations
which induce an even permutation on the set of directions.

Let $\alpha \in G$ be a direction preserving collineation which
induces $\rho$ on $G$ and put $r=\alpha^{-1} \rho$. Obviously,
$\Gamma^+=\langle G, r \rangle$, hence $r \in Z(\Gamma^+)$.
Moreover, since $\Gamma^+$ is invariant under $\sigma$, we have
$r^\sigma \in Z(\Gamma^+)$.

Let $\tau$ be a Bol reflection whose axis is parallel to the axis
of $\sigma$. On the one hand, $\sigma \tau \in G$ and
\[ \sigma \tau = r^{-1} \sigma \tau r =  (\sigma^r \sigma)(\sigma \tau^r)\]
holds. On the other hand, $\sigma^r \sigma = r^{-1} r^\sigma \in
Z(\Gamma^+)$ and $\sigma \tau^r \in \Gamma^+$. Therefore,
\[ (\sigma \tau)^3 = (\sigma^r \sigma)^3 \, (\sigma \tau^r)^3 = \id \]
by the modified triality property, cf. Lemma \ref{l322}.

Assume now that the axis of the Bol reflections $\sigma$ and
$\tau$ are vertical with equation $X=e$ and $X=a$. As we have seen
in Section \ref{subsec:bolreflcoord}, the coordinate forms of
these maps are $(x,y)\sigma = (x^{-1}, xy)$ and $(x,y)\tau =
(ax^{-1}a, a^{-1}(xy))$. This implies $(x,y) \sigma \tau =
(axa,a^{-1}y)$ and $(x,y)(\sigma\tau)^3 = (a^3 x a^3, a^{-3}y)$.
By $(\sigma\tau)^3=\id$, we have $a^3=1$. Since we chose $\tau$
arbitrarily, $L$ must be of exponent $3$.
\end{proof}
\begin{corollary}
If $G$ is a finite group with triality which determines a
noncommutative simple Moufang loop then all triality automorphisms
are outer.
\end{corollary}
\begin{proof}
Assume that $\sigma$ is an inner automorphism. Then so are
$\sigma^\rho$ and $\rho = \sigma^\rho \sigma$. We have the same
implication if we suppose $\sigma\rho$ or $\rho\sigma$ to be
inner. In any case, $\rho$ will be inner and $L$ will be a finite
Moufang loop of exponent $3$. By [13, Thm. 4], $L$ is either not
simple or commutative.
\end{proof}
\begin{theorem}[Liebeck's Theorem \cite{Liebeck}]
The only finite simple groups with triality are the simple groups
$(P\Omega_8^+(q),S)$. The triality automorphisms are uniquely
determined up to conjugation. $($They are the so called graph
automorphisms of $P\Omega_8^+(q).)$
\end{theorem}
\begin{corollary}[Thm. \cite{Liebeck}]
The only nonassociative finite simple Moufang loops are the Paige
loops $M^*(q)=M^*(GF(q))$, where $q$ is a prime power.
\end{corollary}

\section{Automorphism groups of Paige loops over perfect fields}

Now that we have found all nonassociative finite simple Moufang
loops, we will determine their automorphism groups. In fact, we
will determine $\Aut{M^*(F)}$ whenever $F$ is perfect. Recall that
a field of characteristic $p$ is \emph{perfect} if the Frobenius
map $x\mapsto x^p$ is an automorphism of $F$. The approach here is
based on \cite{NagyVojtC}.

\subsection{The automorphisms of the split octonion algebras}

Let $C$ be a composition algebra over $F$. A map $\alpha:C\to C$
is a \emph{linear automorphism} (resp.\ \emph{semilinear
automorphism}) of $C$ if it is a bijective $F$-linear (resp.\
$F$-semilinear) map preserving the multiplication, i.e.,
satisfying $(uv)\alpha=(u\alpha)(v\alpha)$ for every $u$, $v\in
C$. It is well known that the group of linear automorphisms of
$\octo{F}$ is isomorphic to the Chevalley group $G_2(F)$, cf.\
[11, Sec.\ 3], \cite[Ch.\ 2]{SpringerVeldkamp}. The group of
semilinear automorphisms of $\octo{F}$ is therefore isomorphic to
$G_2(F)\rtimes\Aut{F}$.

Since every linear automorphism of a composition algebra is an
isometry \cite[Sec.\ 1.7]{SpringerVeldkamp}, it induces an
automorphism on the loops $M(F)$ and $M^*(F)$. The following
result---that is interesting in its own right---shows that every
element of $\octo{F}$ is a sum of two elements of norm one.
Consequently, $\Aut{\octo{F}}\le\Aut{M^*(F)}$.

\begin{theorem}[Thm.\ 3.3 \cite{Vojtech}]\label{Th:SumsOfTwo}
Let $F$ be any field and $\octo{F}$ the split octonion algebra
over $F$. Then every element of $\octo{F}$ is a sum of two
elements of norm one.

\begin{proof}
As before, we identify $\octo{F}$ with the Zorn vector matrix
algebra over $F$, where the norm is given by the determinant. Let
\begin{displaymath}
    x=\vm{a}{\alpha}{\beta}{b}
\end{displaymath}
be an element of $\octo{F}$. First assume that $\beta\ne 0$. Note
that for every $\lambda\in F$ there is $\gamma\in F^3$ such that
$\gamma\cdot\beta=\lambda$. Pick $\gamma\in F^3$ so that
$\gamma\cdot\beta=a+b-ab+\alpha\cdot\beta$. Then choose
$\delta\in\gamma^\perp\cap\alpha^\perp\ne 0$. This choice
guarantees that $(a-1)(b-1)-(\alpha-\gamma)\cdot(\beta-\delta)
=ab-a-b+1-\alpha\cdot\beta+\gamma\cdot\beta=1$. Thus
\begin{displaymath}
    \vm{a}{\alpha}{\beta}{b}=\vm{1}{\gamma}{\delta}{1}
        +\vm{a-1}{\alpha-\gamma}{\beta-\delta}{b-1}
\end{displaymath}
is the desired decomposition of $x$ into a sum of two elements of
norm $1$. Note that the above procedure works for every $\alpha$.

Now assume that $\beta=0$. If $\alpha\ne 0$, we use a symmetrical
argument as before to decompose $x$. It remains to discuss the
case when $\alpha=\beta=0$. Then the equality
\begin{displaymath}
    \vm{a}{0}{0}{b}=\vm{a}{(1,0,0)}{(-1,0,0)}{0}+\vm{0}{(-1,0,0)}{(1,0,0)}{b}
\end{displaymath}
does the job.
\end{proof}
\end{theorem}

An automorphism $f\in\Aut{M^*(F)}$ will be called
\emph{$($semi$)$linear} if it is induced by a (semi)linear
automorphism of $\octo{F}$.

\subsection{Geometric description of loop automorphisms}

By considering extensions of automorphisms of $M^*(F)$, it was
proved in \cite{Vojtech} that $\Aut{M^*(2)}$ is isomorphic to
$G_2(2)$. The aim of this section is to generalize this result
(although using different techniques) and prove that every
automorphism of $\Aut{M^*(F)}$ is semilinear, provided $F$ is
perfect. We reach this aim by identifying $\Aut{M^*(F)}$ with a
certain subgroup of the automorphism group of the group with
triality associated with $M^*(F)$.

To begin with, we recall the geometrical characterization of
automorphisms of a loop, as promised in Subsection \ref{Ss:Coli}.

\begin{lemma}[Thm.\ 10.2 \cite{BarlStram}]
\label{Lm:GeomChar} Let $L$ be a loop and $\mathcal N$ its
associated $3$-net. Any direction preserving collineation which
fixes the origin of $\mathcal N$ is of the form $(x,\,y)\mapsto
(x\alpha,\,y\alpha)$ for some $\alpha\in\Aut{L}$. Conversely, the
map $\alpha:L\to L$ is an automorphism of $L$ if and only if
$(x,\,y)\mapsto (x\alpha,\,y\alpha)$ is a direction preserving
collineation of $\mathcal N$.
\end{lemma}

We will denote the map $(x,\,y)\mapsto (x\alpha,\,y\alpha)$ by
$\varphi_\alpha$. Before reading any further, recall Propositions
\ref{p324} and \ref{p325}.

\begin{proposition}\label{Pr:Centralizer}
Let $L$ be a Moufang loop and $\mathcal N$ its associated $3$-net.
Let $M$ be the group of collineations generated by the Bol
reflections of $\mathcal N$, $M_0$ the direction preserving part
of $M$, and $S\cong S_3$ the group generated by the Bol
reflections whose axis contains the origin of $\mathcal N$. Then
$\Aut{L}\cong C_{\Aut{M_0}}(S)$.
\end{proposition}

\begin{proof}
As the set of Bol reflections of $\mathcal N$ is invariant under
conjugations by colli\-neations, every element
$\varphi\in\Coll{\mathcal N}$ normalizes the group $M_0$ and
induces an automorphism $\widehat{\varphi}$ of $M$. It is not
difficult to see that $\varphi$ fixes the three lines through the
origin of $\mathcal N$ if and only if $\widehat{\varphi}$
centralizes (the involutions of) $S$.

Pick $\alpha\in\Aut{L}$, and let $\widehat{\varphi_\alpha}$ be the
automorphism of $M_0$ induced by the collineation
$\varphi_\alpha$. As $\varphi_\alpha$ fixes the three lines
through the origin, $\widehat{\varphi_\alpha}$ belongs to
$C_{\Aut{M_0}}(S)$, by the first paragraph.

Conversely, an element $\psi\in C_{\Aut{M_0}}(S)$ normalizes the
conjugacy class of $\sigma$ in $M_0 S$ and preserves the incidence
structure defined by the embedding of $\mathcal N$. This means
that $\psi=\widehat{\varphi}$ for some collineation
$\varphi\in\Coll{\mathcal N}$. Now, $\psi$ centralizes $S$,
therefore $\varphi$ fixes the three lines through the origin. Thus
$\varphi$ must be direction preserving, and there is
$\alpha\in\Aut{L}$ such that $\varphi=\varphi_\alpha$, by Lemma
\ref{Lm:GeomChar}.
\end{proof}

\subsection{The automorphisms of Paige loops}

\begin{theorem}\label{Th:Main}
Let $F$ be a perfect field. Then the automorphism group of the
nonassociative simple Moufang loop $M^*(F)$ constructed over $F$
is isomorphic to the semidirect product $G_2(F)\rtimes\Aut{F}$.
Every automorphism of $M^*(F)$ is induced by a semilinear
automorphism of the split octonion algebra $\octo{F}$.
\end{theorem}
\begin{proof}
We fix a perfect field $F$, and assume that all simple Moufang
loops and Lie groups mentioned below are constructed over $F$.

The group with triality associated with $M^*$ is the
multiplicative group $\Mlt{M^*}\cong D_4$, and the graph
automorphisms of $D_4$ are exactly the triality automorphisms of
$M^*$ (cf. \cite{Freudenthal}, \cite{Doro}). To be more precise,
Freudenthal proved this for the reals and Doro for finite fields,
however they based their arguments only on the root system and
parabolic subgroups, and that is why their result is valid over
any field.

By \cite{Freudenthal}, $C_{D_4}(\sigma) = B_3$, and by
\cite[Lemmas 4.9, 4.10 and 4.3\,]{Liebeck}, $\;C_{D_4}(\rho) =
G_2$. As $G_2<B_3$, by \cite[p.\ 28]{Gorenstein}, we have
$C_{D_4}(S_3) = G_2$.

Since $F$ is perfect, $\Aut{D_4}$ is isomorphic to $\Delta \rtimes
(\Aut{F} \times S_3)$, by a result of Steinberg (cf.\
\cite[Chapter 12]{Carter}). Here, $\Delta$ is the group of the
inner and diagonal automorphisms of $D_4$, and $S_3$ is the group
of graph automorphisms of $D_4$. When $\chr{F}=2$ then no diagonal
automorphisms exist, and $\Delta=\Inn{D_4}$. When $\chr{F}\neq 2$
then $S_3$ acts faithfully on $\Delta / \Inn{D_4} \cong C_2 \times
C_2$. Hence, in any case, $C_{\Delta}(S_3) = C_{D_4}(S_3)$.
Moreover, for the field and graph automorphisms commute, we have
$C_{\Aut{D_4}}(S_3) = C_{D_4}(S_3) \rtimes \Aut{F}$.

We have proved $\Aut{M^*}\cong G_2\rtimes \Aut{F}$. The last
statement follows from the fact that the group of linear
automorphisms of the split octonion algebra is isomorphic to
$G_2$.
\end{proof}

\section{Related results, prospects and open problems}

We conclude with a few results and open problems concerning simple
Moufang loops.

\subsection{Generators for finite Paige loops}

It is well known that every finite simple group is generated by at
most $2$ elements. This result requires the classification of
finite simple groups, and was finalized in \cite{AschGur}. Since
any two elements of a Moufang loop generate a subgroup, no
nonassociative Moufang loop can be $2$-generated. The following
theorem can be proved using some classical results on generators
of groups $SL(2,q)$, cf.\ \cite{PJGT}:

\begin{theorem}\label{Th:Gens}
Every Paige loop $M^*(q)$ is $3$-generated. When $q>2$, the
generators can be chosen as
\begin{displaymath}
    \vm{0}{e_1}{-e_1}{\lambda},
    \quad
    \vm{0}{e_2}{-e_2}{\lambda},
    \quad
    \vm{\lambda}{0}{0}{\lambda^{-1}},
\end{displaymath}
where $\lambda$ is a primitive element of $GF(q)$, and $e_i$ is
the $3$-vector whose only nonzero coordinate is in position $i$
and is equal to $1$. When $q=2$, the generators can be chosen as
\begin{displaymath}
    \vm{1}{e_1}{e_1}{0},\quad
    \vm{1}{e_2}{e_2}{0},\quad
    \vm{0}{e_3}{e_3}{1}.
\end{displaymath}
\end{theorem}

\subsection{Generators for integral Cayley numbers of norm one}

Let $C=(C,N)$ be a real composition algebra. The \emph{set of
integral elements} of $C$ is the maximal subset of $C$ containing
$1$, closed under multiplication and subtraction, and such that
both $aN$ and $a+\overline{a}$ are integers for each $a$ in the
set.

Let $\mathbb R$, $\mathbb C$, $\mathbb H$, $\mathbb O$ be the
classical real composition algebras, i.e., those obtained from
$\mathbb R$ by the Cayley-Dickson process with parameter
$\lambda=-1$. The real octonions $\mathbb O$ are often called
\emph{Cayley numbers}. For $C\in\{\mathbb R$, $\mathbb C$,
$\mathbb H\}$, there is a unique set of integral elements of $C$.
(For instance, when $C=\mathbb C$ this set is known as the
Gaussian integers.) When $C=\mathbb O$, there are seven such sets,
all isomorphic, as Coxeter showed in \cite{Coxeter}.

We use the notation of \cite{Coxeter} here. Let $J$ be one of the
sets of integral elements in $\mathbb O$, and let $J'=\{x\in J|\
xN=1\}$. Then $|J'|=240$, and $J'/\{1,\,-1\}$ is isomorphic to
$M^*(2)$. (This may seem strange, however, $M^*(2)$ is a subloop
of any $M^*(q)$, too.) Hence, by Theorem \ref{Th:Gens},
$J'/\{1,\,-1\}$ must be $3$-generated. Let $i$, $j$, $k$ be the
usual units in $\mathbb H$, and let $e$ be the unit that is added
to $\mathbb H$ when constructing $\mathbb O$. Following Dickson
and Coxeter, let $h=(i+j+k+e)/2$. Then one can show that $i$, $j$
and $h$ generate $J'/\{1,\,-1\}$ (multiplicatively). Since
$i^2=-1$, it follows that every set of integral elements of unit
norm in $\mathbb O$ is $3$-generated, too. See \cite{PJA} for
details.

\subsection{Problems and Conjectures}

\begin{problem} Find a presentation for $M^*(q)$ in the variety of Moufang
loops, possibly based on the generators of Theorem
$\ref{Th:Gens}$.
\end{problem}

\begin{problem} Find $($necessarily infinite$)$ nonassociative simple Moufang
loops that are not Paige loops.
\end{problem}
\begin{conjecture} Let $L$ be a nonassociative simple Moufang loop and let
$H=\mathrm{Mlt}(L)_e$ be the stabilizer of the neutral element in
the multiplication group of $L$. Then $H$ is simple.
\end{conjecture}
\begin{problem} \label{prob:fn} Find a function $f: \mathbb  N \to \mathbb
N$ such that the order of the multiplication group of a Moufang
loop of order $n$ is less that $f(n)$.
\end{problem}
For the finite Paige loop $M^*(q)$, we have
\begin{eqnarray*}
    |M^*(q)|&=&\frac{1}{d} q^3(q^4-1),\\
    |P\Omega^+_8(q)| &=& \frac{1}{d^2} q^{12} (q^2-1) (q^4-1)^2 (q^6-1),
\end{eqnarray*}
where $d=(2,q-1)$. Hence $|\Mlt(M^*(q))| < 4 |M^*(q)|^4$ holds.
This motivated us to state:
\begin{conjecture}
The function $f(n)=4n^4$ solves Problem $\ref{prob:fn}$.
\end{conjecture}

\small
\bibliographystyle{plain}

\noindent \footnotesize{\hfill Received \ April 29, 2003

\medskip
\begin{tabular}{lll}
        G\'abor P. Nagy&\quad&Petr Vojt\v{e}chovsk\'y\\
        Bolyai Institute&\quad&Department of Mathematics\\
        University of Szeged&\quad&University of Denver\\
        Aradi v\'ertan\'uk tere 1&\quad&2360 S Gaylord St\\
        H-6720 Sze\-ged&\quad&Denver, Colorado 80208\\
        Hungary&\quad&U.S.A.\\
        e-mail: nagyg@math.u-szeged.hu&\quad&e-mail: petr@math.du.edu
\end{tabular}
}
\end{document}